\documentclass[oneside,12pt]{amsart}
\usepackage{math,wrapfig}
\usepackage[all]{xy}
\newcommand\8{\mathsf{0}}
\newcommand\9{\mathsf{1}}

\newcommand\commpar[1]{\hbox{\vline~\parbox[t]{\marginparwidth}{\small\sf #1}}}
\newcommand\comment[1]{\begin{wrapfigure}{r}{\marginparwidth}%
    \vspace{-2.5ex}\commpar{#1}\end{wrapfigure}}
\newcommand\ncomment[2]{\begin{wrapfigure}[#1]{r}{\marginparwidth}%
    \vspace{-2.5ex}\commpar{#2}\end{wrapfigure}}

\begin{document}
\title[The dimension series of the BSV group]{The
  $2$-dimension series of the just-nonsolvable BSV group}
\author{Laurent Bartholdi}
\date{February 23, 2001; typeset \today}
\email{laurent@math.berkeley.edu}
\address{11 ch. de la Barillette, 1260 Nyon, Switzerland}
\thanks{The author thanks the ``Swiss National Foundation for
  Scientific Research''}
\dedicatory{To Said N. Sidki, born January 23, 1941 in Al Quds
  (Palestine), for his 60th birthday}
\keywords{Lie algebra; Lower Central Series; Dimension Series}
\subjclass{\textbf{20F14} (Derived series, central series, and
  generalizations); \textbf{20F40} (Associated Lie structures);
  \textbf{17B70} (Graded Lie (super)algebras); \textbf{20E08} (Groups
  acting on trees)}
\begin{abstract}
  \hspace{-1cm}\parbox[t]{75mm}{\hspace{1cm}I compute the structure of
    the restricted $2$-algebra associated to a group first described
    by Andrew Brunner, Said Sidki and Ana Cristina Vieira, acting on
    the binary rooted tree~\cite{brunner-s-v:nonsolvable}. I show that
    its width is unbounded, growing logarithmically, and obeys a
    simple rule. As a consequence, the dimension of
    $\varpi^n/\varpi^{n+1}$ (where $\varpi<\F_2\Gamma$ is the
    augmentation ideal), is $p(0)+\dots+p(n)$, the total number of
    partitions of numbers up to $n$.}\quad\parbox[t]{35mm}{\tiny This
    paper contains many marginal notes.  None, including this one, is
    necessary for the understanding of the results.  However I hope
    that they will provide some insight into the motivation and
    process of discovery, as opposed to a description of the end
    result.}
\end{abstract}
\maketitle

\section{Introduction}
\ncomment4{The standard Lie algebra is much more complicated;
  considering this restricted Lie algebra suppresses much of the
  $2$-torsion.}  The purpose of this paper is to construct the
restricted Lie algebra associated to the Brunner-Sidki-Vieira group
$\Gamma$. I decide to state the auxiliary results in as little
generality as possible, aiming to arrive as quickly as possible at an
explicit description.

\subsection{The group}
\comment{$\Gamma$ is universal in the following sense: any
  finite-state state-closed infinite cyclic group of the automorphism
  group of the binary tree is generated by an $n$-th root of $\tau$ or
  $\mu$ for odd $n$~\cite{nekrashevych-s:sc}. This singles out $\tau$
  and $\mu$ as ``special'' generators.  $\tau$ is called the ``dyadic
  adding machine'', or ``odometer''.} Let $\Sigma=\{0,1\}$ be a
two-letter alphabet.  The \emph{rooted binary tree} is the free monoid
$\Sigma^*$.  The group $\Gamma$ acts on $\Sigma^*$, and is generated
by two elements $\tau,\mu$ defined as follows:

\begin{xalignat*}{2}
  \emptyset^\tau&=\emptyset, & \emptyset^\mu&=\emptyset;\\
  (0\sigma_2\dots\sigma_n)^\tau&=1\sigma_2\dots\sigma_n, &
  (0\sigma_2\dots\sigma_n)^\mu &=1\sigma_2\dots\sigma_n;\\
  (1\sigma_2\dots\sigma_n)^\tau&=0(\sigma_2\dots\sigma_n)^\tau, &
  (1\sigma_2\dots\sigma_n)^\mu &=0(\sigma_2\dots\sigma_n)^{\mu^{-1}}.\\
\end{xalignat*}

\ncomment3{The growth of $\Gamma$ is unknown, and that motivated my
  computations --- see equation~(\ref{eq:growth}). It is known that
  $\Gamma$ doesn't contain non-abelian free
  subgroups~\cite{silva:phd}, but it is unknown whether it contains
  non-abelian free monoids.}  It is known that this group is infinite,
just-nonsolvable, and torsion-free~\cite{brunner-s-v:nonsolvable}.

\vspace{3\baselineskip}

\ncomment3{also called ``Brauer'', ``Jennings'', ``Lazard'', and
  ``Zassenhaus'' series; $2$ can be replaced by any prime.}
\subsection{The dimension series} The \emph{$2$-dimension series} of
a group $\Gamma$ is defined as follows: $\Gamma_1=\Gamma$, and
$\Gamma_n=[\Gamma,\Gamma_{n-1}]\Gamma_{\lceil n/2 \rceil}^2$, where
$H^2$ denotes the subgroup of $H$ generated by its squares. It can
alternately be described, by a result of Lazard~\cite{lazard:nilp}, as
\[\Gamma_n = \prod_{2^ji\ge n}\gamma_i(\Gamma)^{2^j},\]

\comment{it was thought at first that this was an alternate definition
  of the lower central series, for $2$-groups. The correspondence
  between the $\Gamma_i$ and the $\gamma_i(\Gamma)$ is a deep
  question.}
\noindent $\gamma_i(\Gamma)$ being the lower central series, or as
\[\hspace{-65mm}\Gamma_n = \{g\in\Gamma|\,g-1\in\varpi^n\},\]
where $\varpi$ is the augmentation ideal of the group algebra
$\F_2\Gamma$.

\comment{$\Gamma_n/\Gamma_{n+1}$ is an elementary abelian $2$-group,
  by construction.}
\subsection{The Lie algebra} The \emph{$2$-Lie algebra} of $\Gamma$
is the restricted graded Lie algebra
\[\Lie(\Gamma)=\bigoplus_{n\ge1}\Gamma_n/\Gamma_{n+1}.\]
Its Lie bracket is induced by commutation in $\Gamma$, and its
Frobenius endomorphism is induced by squaring in $\Gamma$.

The \emph{degree} of $g\in\Gamma$ is the maximal
$n\in\N\cup\{\infty\}$ such that $g$ belongs to $\Gamma_n$.
A basis of the rank-$n$ module of $\Lie(\Gamma)$ can be found among
elements of degree $n$; and since $\Gamma$ is residually-$2$, the only
element of infinite degree is $1$.

\comment{However, the width could also be defined with respect to the
  lower central series; I don't know whether these definitions are
  equivalent.}  Following~\cite{klass-lg-p:fw}, say $\Gamma$ has
\emph{finite width} if there is a constant $W$ such that
$\ell_n:=\dim\Gamma_n/\Gamma_{n+1}\le W$ holds for all $n$.

Consider the graded algebra
\[\overline{\F_2\Gamma}=\bigoplus_{n=0}^\infty\varpi^n/\varpi^{n+1}.\]
A fundamental result of Daniel Quillen~\cite{quillen:ab} implies that
$\overline{\F_2\Gamma}$ is the enveloping $2$-algebra of $\Lie(\Gamma)$.

\comment{Of course, Jennings did that before Quillen, and quite
  unaware of the fact that his construction corresponded to PBW!}
Using the Poincar\'e-Bikhoff-Witt isomorphism, Stephen
Jennings~\cite{jennings:gpring} then showed that
\vspace{1ex}
\begin{equation}\label{eq:jennings}
  \sum_{n\ge0}\dim(\varpi^n/\varpi^{n+1})\hbar^n
  =\prod_{n\ge1}(1+\hbar^n)^{\ell_n}.
\end{equation}

\comment{This shows that a residually-$2$ group either has polynomial
  growth (if $\ell_n=0$ for some $n$) or growth at least $e^{\sqrt n}$
  (if $\ell_n\ge1$ for all $n$).} It is also
known~\cite[Lemma~2.5]{bartholdi-g:lie} that if $f_n$ denotes the
number of group elements of $\Gamma$ of length at most $n$ in the
generators, then
\begin{equation}\label{eq:growth}
  f_n\ge\dim\F_2\Gamma/\varpi^{n+1}=\sum_{i=0}^n\dim(\varpi^n/\varpi^{n+1}).
\end{equation}

\comment{At present, only the first four successive quotients of the
  lower central series are known: $\gamma_1/\gamma_2=\Z^2$,
  $\gamma_2/\gamma_3=\Z$ and
  $\gamma_3/\gamma_4=\gamma_4/\gamma_5=\Z/8\Z$.}
\subsection{The main result}
Along with a description of $\Lie(\Gamma)$, I shall show:
\begin{thm}\label{thm:main}
  The ranks of successive quotients,
  $\ell_n=\dim_{\F_2}\Gamma_n/\Gamma_{n+1}$, satisfy
  \[\ell_n=\begin{cases}i+2& \text{ if }n=2^i\text{ for some }i;\\
    \max\setsuch{i+1}{2^i\text{ divides }n}& \text{ otherwise}.
  \end{cases}\]
  As a consequence,
  $\dim\varpi^n/\varpi^{n+1}=\sum_{i=0}^np(n)$, the number of
  partitions of numbers up to $n$, and the word growth of
  $\Gamma$ is at least $e^{\sqrt n}$.
\end{thm}

\section{Lie graphs and Branch Portraits}
\ncomment2{This is in fact just a simple representation of the
  structure constants of $\Lie(\Gamma)$.}  The structure of
$\Lie(\Gamma)$ can best be described using Lie graphs, introduced
in~\cite{bartholdi:lcs}. The \emph{Lie graph} of $\Lie(\Gamma)$ has as
vertex set a basis $V$ of $\Lie(\Gamma)$; in the favourable case that
$[v,\tau]$, $[v,\mu]$ and $v^2$ are basis elements for all $v\in
V$,\par
\ncomment2{This is the first group for which this miracle happens!}
\noindent the Lie graph has for all $v\in V$ an arrow labelled $\tau$
from $v$ to $[v,\tau]$, one labelled $\mu$ leading to $[v,\mu]$ and
one labelled $2$ leading to $v^2$.

\ncomment3{This property imposes a great amount of rigidity on
  $\Gamma$; such groups are usually referred to as ``weakly
  branch''~\cite{bartholdi-g:spectrum}.}
\subsection{Branch portraits} The important ``branching property'' of
$\Gamma$ is that $\Gamma'$ contains $\Gamma'\times\Gamma'$, where in
this last group the left and right factors act on $0\Sigma^*$ and
$1\Sigma^*$ respectively; and furthermore $\Gamma'=\langle
c\rangle(\Gamma'\times\Gamma')$, with $c=[\tau,\mu]$.
\newpage
\begin{wrapfigure}[4]{r}{\marginparwidth}
  \commpar{The advantage of this labelling is that on prescribing any
    finite set of label values one still defines elements of $\Gamma$.
    This property is not shared by the labelling below.}
\end{wrapfigure}
Therefore every $g\in\Gamma$ can be represented as
$g=\tau^i\mu^jg_\emptyset$ for some $i,j\in\Z$ and
$g_\emptyset\in\Gamma'$; then inductively each $g_\sigma\in\Gamma'$
can be represented as
$g_\sigma=c^{k_\sigma}(g_{0\sigma},g_{1\sigma})$, and so every
$g\in\Gamma$ gives rise to a labelling of the tree $\Sigma^*$, with
label $(i,j,k_\emptyset)\in\Z^3$ at the root vertex $\emptyset$ and
label $k_\sigma\in\Z$ at the vertex $\sigma$.

Consider the group, written $\overline\Gamma$, of all infinite
expressions of the form
\[\tau^*\mu^*c^*(c^*,c^*)(c^*,c^*,c^*,c^*)\cdots\]
for any choices of $*\in\Z$; this is an uncountable group containing
$\Gamma$. The dimension series of $\overline\Gamma$ coincides with
that of $\Gamma$, so we may perform the computations in
$\overline\Gamma$, written again $\Gamma$ from now on.

\ncomment3{This gives a labelling of the tree with labels in $\Z/2$ ---
  but not all labels define an element of $\Gamma$, or even of its
  completion.}  Elements of $\Gamma$ may also be decomposed according
to their action on $\Sigma^*$; to wit, given $g\in\Gamma$, choose
$\epsilon^i\in\mathsf{Sym}(2)=\langle\epsilon\rangle$ such that
$g/\epsilon^i$ fixes the first letter of all words in $\Sigma^*$; then
$g/\epsilon^i$ acts both on $0\Sigma^*$ and $1\Sigma^*$ as elements of
$\Gamma$, written $g_0$ and $g_1$. In this way every $g\in\Gamma$ has
an expression of the form
\[\epsilon^*(\epsilon^*,\epsilon^*)
(\epsilon^*,\epsilon^*,\epsilon^*,\epsilon^*)\cdots;\]
I shall implicitly write elements in either form $g$ or
$(g_0,g_1)\epsilon^i$. Using that notation,
\begin{align*}
  \tau&=(1,\tau)\epsilon,\\
  \mu&=(1,\mu^{-1})\epsilon,\\
  c&=\epsilon(1,\tau^{-1})\epsilon(1,\mu)(1,\tau)\epsilon
  (1,\mu^{-1})\epsilon=((\mu\tau)^{-1},\mu\tau).
\end{align*}

\comment{These maps will be used to create arbitrary branch portraits;
  it's much more clever to use these elements $(g,1)$ and $(g,g)$
  rather than the ``obvious'' $(g,1)$ and $(1,g)$.}
\subsection{The basis} Recalling that $\Gamma'$ contains
$\Gamma'\times\Gamma'$, define the following two endomorphisms of
$\Gamma'$:
\[\hspace{-60mm}\8(g)=(g,1),\qquad\9(g)=(g,g).\]

\vspace{1ex}

\begin{thm}
  The following elements form a basis of $\Lie(\Gamma)$:
  \begin{itemize}
  \item $\tau^{2^n}$ for all $n\in\N$, of degree $2^n$;
  \item $\mu^{2^n}$ for all $n\in\N$, of degree $2^n$;
  \item $W(c)^{2^n}$ for all $W\in\{\8,\9\}^m$ and $n\in\N$, of
    degree
    \[\hspace{16mm} 2^n\big(\sum_{i=1}^m
    W_i2^{i-1}+2^m+1\big).\hspace{16mm}\commpar{That's really reading
      `$\9W\8\dots\8\9$' in base $2$.}\]
  \end{itemize}

  The arrows in the Lie graph are
  \[\xymatrix@R-5mm{{\tau^{2^n}}\ar[r]^\mu & {\8^n(c)} &
    {\mu^{2^n}}\ar[r]^\tau & {\8^n(c)}\\
    & {\9^n(c)}\ar[r]^{\tau,\mu} & {\8^{n+1}(c)}\\
    {v^{2^n}}\ar[r]^2 & {v^{2^{n+1}}\makebox[0mm][l]{ for
        all basis vectors $v$}}\\
    {\9^n\8W(c)}\ar[r]^{\tau,\mu} & {\8^n\9W(c)\makebox[0mm][l]{
        for all $W\in\{\8,\9\}^*$}}}\]
\end{thm}

\begin{figure}
  {\tiny\xymatrix@d{
    *+[o][F-]{1} & *+[o][F-]{2} & *+[o][F-]{3} & *+[o][F-]{4} & *+[o][F-]{5} & 
    *+[o][F-]{6} & *+[o][F-]{7} & *+[o][F-]{8} & *+[o][F-]{9} & *+[o][F-]{10}\\
    & & & & & & & {c^4}\ar@{.}[rrr] & & & {}\\
    & & & {c^2}\ar[urrrr]^2 & & {\8(c)^2} & & {\9(c)^2} & &
    {\8\8(c)^2}\ar@{.}[r] & {}\\
    & {c}\ar[urr]^2\ar[r]_{\tau,\mu} &
    {\8(c)}\ar[r]_{\tau,\mu}\ar[urrr]^2 &
    {\9(c)}\ar[r]_{\tau,\mu}\ar[urrrr]^2 &
    {\8\8(c)}\ar[r]_{\tau,\mu}\ar[urrrrr]^2 &
    {\9\8(c)}\ar[r]_{\tau,\mu} &
    {\8\9(c)}\ar[r]_{\tau,\mu} &
    {\9\9(c)}\ar[r]_{\tau,\mu} &
    {\8\8\8(c)}\ar[r]_{\tau,\mu} &
    {\9\8\8(c)}\ar@{.}[r] & {}\\
    {\tau}\ar[ur]^{\mu}\ar[r]^2 &
    {\tau^2}\ar[ur]^{\mu}\ar[rr]^2 & &
    {\tau^4}\ar[ur]^{\mu}\ar[rrrr]^2 & & & &
    {\tau^8}\ar[ur]^{\mu}\ar@{.}[rrr] & & & {}\\
    {\mu}\ar[uur]_(0.4){\tau} \ar[r]_2 &
    {\mu^2}\ar[uur]_(0.4){\tau} \ar[rr]_2 & &
    {\mu^4} \ar[uur]_(0.4){\tau} \ar[rrrr]^2 & & & &
    {\mu^8} \ar[uur]_(0.4){\tau} \ar@{.}[rrr] & & & {}}}
  \caption{The beginning of the Lie graph of $\Lie(\Gamma)$.}
\end{figure}
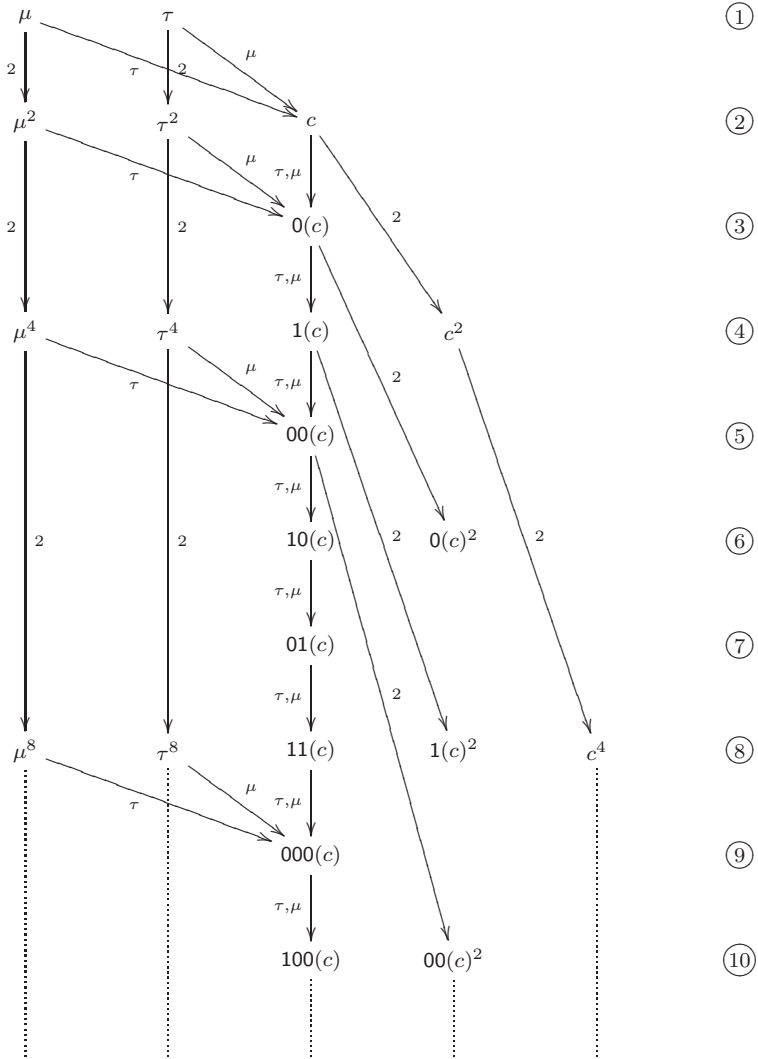

\newpage
{\itshape \proofname.}  The strategy of the proof
strategy is simple, and follows ideas from~\cite{bartholdi:lcs}.

\ncomment4{These were guessed using \textsf{Gap}~\cite{gap:manual}, by
  computing the $2$-series of the quotient of $\Gamma$ acting on the
  eigth level $\Sigma^8$ of the tree. The use of $\8$ and $\9$ is
  motivated by the fact that $[\8v,\tau]\equiv1v$.}  First the degrees
of the basis vectors are as claimed: for this it suffices to check
that the commutation relations are as described.  Compute:

\vspace{1ex}

\begin{align*}
  [\mu^{2^n},\tau]&=(\mu^{2^{n-1}},\mu^{2^{n-1}})\epsilon
  (1,\tau^{-1})(\mu^{-2^{n-1}},\mu^{-2^{n-1}})(1,\tau)\epsilon\\
  &=([\mu^{-2^{n-1}},\tau],1)=\8[\mu^{-2^{n-1}},\tau]=\dots
  =\8^n(c)^{(-1)^n}\equiv\8^n(c);\\
  [\9^n\8W(c),\tau]&=(\9^{n-1}\8W(c)^{-1},\9^{n-1}\8W(c)^{-1})
  \epsilon(1,\tau^{-1})\\
  &\hspace{10em} (\9^{n-1}\8W(c),\9^{n-1}\8W(c))(1,\tau)\epsilon\\
  &=([\9^{n-1}\8W(c),\tau],1)=\dots=\8^n[\8W(c),\tau]\\
  &=\8^n(W(c)^{-1},1)\epsilon(1,\tau^{-1})(W(c),1)(1,\tau)\epsilon\\
  &=\8^n\9W(c)/\8^{n+1}W(c)^2\equiv\8^n\9W(c);\\
  [\9^n(c),\tau]&=(\9^{n-1}(c)^{-1},\9^{n-1}(c)^{-1})\epsilon
  (1,\tau^{-1})(\9^{n-1}(c),\9^{n-1}(c))(1,\tau)\epsilon\\
  &=([\9^{n-1}(c),\tau],1)=\dots=\8^n[c,\tau]\\
  &=\8^n(\mu\tau,(\mu\tau)^{-1})\epsilon(1,\tau^{-1})
  ((\mu\tau)^{-1},\mu\tau)(1,\tau)\epsilon\\
  &=\8^{n+1}(c^{-\tau^{-1}\mu^{-2}})/\8^n(c)^2\equiv\8^{n+1}(c);
\end{align*}
\ncomment4{It is there that ``dropping'' squares makes the
  computations so much simpler for the dimension series than for the
  lower central series.} and similar results hold for commutation with
$\mu$.  Here the sign `$\equiv$' means congruence modulo terms of
greater degree, i.e.\ commutators or squares of the elements involved.
  
\comment{Imagine the Lie graph as hanging from its top generators
  $\tau,\mu$, with length-$1$ and length-doubling edges attaching the
  basis elements. The depth of an element is the maximal number of
  consecutive edges connecting it to the top. In this image all edges
  are under tension, uniquely defining the Lie graph structure. The
  Lie graph is determined by having degree-compatible edges, and being
  connected. } It remains to check that the vectors of degree $n$ are
independent modulo $\Gamma_{n+1}$; this is done by induction on $n$,
the first two values following from~\cite{brunner-s-v:nonsolvable}.
Consider first odd integers $n>1$, when there is a single element
$W(c)$ to consider.  If $W(c)\in\Gamma_{n+1}$, then since
$W(c)=[V(c),\tau]=[V(c),\mu]$ for some $V(c)$, this would imply
$V(c)\in\Gamma_n$ contrary to induction.
    
Consider now even integers $n=2m$, and the elements $W(c),A^2$ of our
putative basis of degree $n$, for some $A\in\Gamma_m$. Among the
linear combinations in $\Gamma_{n+1}$, If $W(c)\in\Gamma_{n+1}$, then
again this would imply $V(c)\in\Gamma_n$ for some $V(c)$ of lesser
weight; if $A^2\in\Gamma_{n+1}$, then $A\in\Gamma_{m+1}$; if
$W(c)A^2\in\Gamma_{n+1}$, then $V(c)\in\Gamma_n$ and
$A\in\Gamma_{m+1}$. All these conclusions are contrary to induction,
and are the only possible linear combinations among basis elements of
degree $n$.  \qed

The proof of Theorem~\ref{thm:main} now follows; indeed the basis
vectors of degree $2^n$ are $\tau^{2^n}$, $\mu^{2^n}$ and the
$\9^i(c)^{n-i}$ for all $i\in\{0,\dots,n-1\}$, making $n+2$ elements
in all, and the basis vectors of degree $2^nm$, for odd $m$, are the
$\9^iW(c)^{n-i}$ for all $i\in\{0,\dots,n\}$, where $W$ is the unique
word in $\{\8,\9\}^*$ such that $W(c)$ is of degree $m$. Then
by~(\ref{eq:jennings})
\begin{align*}
  \sum_{n\ge0}\dim(\varpi^n/\varpi^{n+1})\hbar^n
  &=\prod_{n\ge0}(1+\hbar^{2^n})\prod_{n\ge0}\prod_{m\ge1}(1+\hbar^{2^nm})\\
  &=\frac1{1-\hbar}\prod_{m\ge1}\frac1{1-\hbar^m}
  =\sum_{0\le i\le n}p(i)\hbar^n.
\end{align*}
\comment{In fact, this lower bound is already known, because $\Gamma$
  is not virtually nilpotent.}  The asymptotics $p(n)\approx
e^{\pi\sqrt{2n/3}}$ are well known to number
theorists~\cite{hardy-r:asymptotic,uspensky:asymptotic}, and imply
by~(\ref{eq:growth}) the lower bound on the growth of $\Gamma$.

\bibliography{mrabbrev,people,math,grigorchuk,bartholdi}
\end{document}